\documentclass[12pt]{amsart}

\usepackage{mathrsfs}
\usepackage{amsmath}
\usepackage{latexsym}
\usepackage{amssymb}
\usepackage{amsfonts}
\usepackage{bm}
\include{PDF}
\usepackage{graphics}
\def\l{\langle} \def\r{\rangle} 
\def\div{\,\,\big|\,\,} 
 \def\ZZ{\mathbb Z}

 \def\val{{\sf val}}

\def\Aut{{\sf Aut}} 

\def\Cos{{\sf Cos}}

\def\HS{{\sf HS}}

\def\Cay{{\sf Cay}} 
\def\D{{\rm D}} 
\def\S{{\rm S}} 
\def\J{{\rm J}} \def\M{{\rm M}}

\def\C{{\rm C}}\def\N{{\rm N}} \def\O{{\bf O}}

\def\Ga{{\it \Gamma}}

\def\a{\alpha}

\def\Mult{{\sf Mult}}

 \def\GL{{\rm GL}}

\def\PSp{{\rm PSp}}

\def\A{{\rm A}}

\def\PSL{{\rm PSL}}\def\PSO{{\rm PSO}}
\def\GL{{\rm GL}} 

 \def\PSU{{\rm PSU}} 
 \def\McL{{\rm McL}} 
 
 \def\F{{\rm F}} \def\D{{\rm D}}

\def\HS{{\rm HS}}

\newtheorem{rmk}{Remark}[section]
\newtheorem{thm}{Theorem}[section]

\newtheorem{problem}[thm]{Problem}

\newtheorem{lemma}{Lemma}[section]

\newtheorem{example}[thm]{Example}
\def\pf{\noindent{\it Proof.} }
\def\qed{\nopagebreak\hfill{\rule{4pt}{7pt}}
\medbreak}

\begin{document}

\title[Symmetric Cayley graphs]{On symmetric Cayley graphs of valency thirteen$^{*}$}
\thanks{2010 MR Subject Classification 20B15, 20B30, 05C25.}
\thanks{$^{*}$
The work was supported by the National Natural Science Foundation of China (11241076, 11861076). }
\thanks{$^\dag$Corresponding author. E-mails: bengong188@163.com (B.G. Lou).}
\author[B.G. Lou, Z. Zuo, B. Ling] {Ben Gong Lou$^{1\dag}$, Zheng Zuo$^{1}$, Bo Ling$^{2}$}

\address{1: School of Mathematics and Statistics\\ Yunnan University,
Kunming, Yunnan 650500, P. R. China}
\email{bengong188@163.com (B.G. Lou)}

\address{2: School of Mathematics and Computer Sciences\\
Yunnan Minzu University\\
Kunming, Yunnan 650504, P.R. China}
\email{bolinggxu@163.com (B. Ling)}


\begin{abstract}
A Cayley graph $\Ga=\Cay(G,S)$ is said to be normal if the right-regular representation of $G$ is normal in $\Aut\Ga$.
In this paper, we investigate the normality problem of the connected 13-valent symmetric Cayley graphs $\Ga$
of finite nonabelian simple groups $G$, where the vertex stabilizer $\A_v$ is soluble for $\A=\Aut\Ga$ and $v\in V\Ga$.
We prove that $\Ga$ is either normal or $G=\A_{12}$, $\A_{38}$, $\A_{116}$, $\A_{207}$, $\A_{311}$, $\A_{935}$ or $\A_{1871}$.
Further, 13-valent symmetric non-normal Cayley graphs of $\A_{38}$, $\A_{116}$ and $\A_{207}$ are constructed.
This provides some more examples of non-normal 13-valent symmetric Cayley graphs of finite nonabelian simple groups
since such graph (of valency 13) was first constructed
by Fang, Ma and Wang in (J. Comb. Theory A 118, 1039--1051, 2011).

\vskip4pt


\end{abstract}

\maketitle
\section{Introduction}
All graphs are assumed to be finite, simple and undirected in this paper.

Let $\Ga$ be a graph. We use $V\Ga$, $E\Ga$
and $\Aut\Ga$ to denote the vertex set, edge set and
automorphism group of $\Ga$, respectively.
Denote $\val\Ga$ the valency of
$\Ga$. Let $X\le\Aut\Ga$. The graph $\Ga$ is said to be {\em $X$-vertex-transitive},
if $X$ is transitive on $V\Ga$. If $X$ is transitive on the set of arcs of $\Ga$,
then $\Ga$ is called an {\em $X$-arc-transitive graph} or an {\em $X$-symmetric graph}.
In particular, if $X=\Aut\Ga$, then $\Ga$ is simply called {\em vertex-transitive} or {\em arc-transitive} (or {\em symmetric}), respectively.

Let $G$ be a finite group with identity $1$, and let $S$ be a subset of $G$ such that
$1\not\in S$ and $S=S^{-1}:=\{x^{-1}\mid x\in S\}$. The  Cayley graph of $G$ with respect to $S$, denoted by $\Cay(G,S)$,
is defined on $G$ such that $g,\,h\in G$ are adjacent if and only if $hg^{-1}\in S$.
For a Cayley graph $\Cay(G,S)$, the underlying group $G$ can be viewed as a regular subgroup of $\Aut\Cay(G,S)$ which acts on $G$ by right multiplication.
Then a Cayley graph $\Ga=\Cay(G,S)$ is said to be {\em normal} if $G$ is normal in $\Aut\Ga$; otherwise, $\Ga$ is called {\em non-normal}.

The concept of normal Cayley graphs was first proposed by M.Y.Xu in \cite{Xu98} and
it plays an important role in determining the full automorphism groups of Cayley graphs.
The Cayley graphs of finite nonabelian simple groups are received most attention in the
literature.
In 1996, C.H.Li \cite{Li96} proved that a connected cubic symmetric Cayley graph of
a nonabelian simple group $G$ is normal except $7$ groups.
On the basis of C.H.Li's result,
S.J.Xu et al. \cite{Xufang,XufangWang} proved that all such graphs are normal except two
Cayley graphs of the alternating group $\A_{47}$.
In 2002, Fang, Praeger and Wang \cite{FPW02} developed a theory for investigating the automorphism
groups of Cayley graphs of nonabelian simple groups, which is then used to characterize
locally primitive Cayley graphs (that is, $(\Aut\Ga)_v$ acts primitively on the neighbourhood
$\Ga(v)$ for a vertex $v$ of $\Ga$) of nonabelian simple groups by \cite{FMW11}.
Further, Fang, Ma and Wang in \cite{FMW11} proved that all but finitely many locally primitive Cayley graphs of valency $d\le20$ or a prime
number of
the finite nonabelian simple groups are normal. Then they proposed the following problem:
\begin{problem}\label{problem1}
Classify non-normal locally primitive Cayley graphs of finite simple groups with valency $d\le20$ or a prime number.
\end{problem}

From the classification of the small valencies, we know that examples of connected symmetric non-normal Cayley graphs
of nonabelian simple groups are very rare (see \cite{DFZ17-2,FangLiXu-01,FPW02,FWZ16,LiLing01} for valency four,
\cite{DFZ17,LL16,Zhou-1} for valency five, \cite{LiJJ-7,PYL18} for valency seven, \cite{LingLou11} for valency eleven).
We concentrate on the 13-valent case in this paper.
The first known example of non-normal 13-valent symmetric Cayley graph of
nonabelian simple group was constructed by Fang, Ma and Wang \cite{FMW11},
that is, the non-normal Cayley graph of $\A_{12}$.
The aim of this paper is to classify the connected non-normal 13-valent symmetric Cayley graphs
with soluble vertex stabilizers on finite nonabelian simple groups.
In particular, we will construct non-normal 13-valent symmetric Cayley graphs on $\A_{38}$, $\A_{116}$ and $\A_{207}$.

Our main result is the following theorem.
\begin{thm}\label{thm1}
Let $G$ be a finite nonabelian simple group, and let $\Ga=\Cay(G,S)$ be a connected 13-valent symmetric Cayley graph of $G$.
Let $\A=\Aut\Ga$ and $\A_{v}$ be the stabilizer of $v$ in $\A$ where $v\in V\Ga$.
If $\A_v$ is soluble, then the following statements hold.
\begin{itemize}
\item[(1)] Either $\Ga$ is a normal Cayley graph or $G=\A_{12}$, $G=\A_{12}$, $\A_{38}$, $\A_{116}$, $\A_{207}$, $\A_{311}$, $\A_{935}$ or $\A_{1871}$. Further,
\item[(2)] there exist connected non-normal 13-valent symmetric Cayley graphs for $G=\A_{12}$, $\A_{38}$, $\A_{116}$ or $\A_{207}$.
\end{itemize}
\end{thm}
\begin{rmk}\rm
(a) The connected non-normal 13-valent symmetric Cayley graph of $\A_{12}$ was constructed
by Fang, Ma and Wang in \cite{FMW11}.

(b) Specific examples of $\A_{38}$, $\A_{116}$ and $\A_{207}$ which satisfy parts (2) are constructed in Section 4.

(c) We do not know whether all connected 13-valent symmetric Cayley graphs of $\A_{312}$, $\A_{936}$ or $\A_{1872}$ are normal.
\end{rmk}

\section{Preliminaries}
We give some necessary preliminary results in this section.

Let $G$ be a group, $g\in G$
and $H$ a subgroup of $G$.
Define the {\it coset graph} $\Cos(G, H, g)$ of $G$ with respect
to $H$ as the graph with vertex set $[G:H]$ (the set of cosets of $H$ in $G$),
and $Hx$ is adjacent to $Hy$ with $x,y\in G$ if and only if
$yx^{-1}\in HgH$. The following lemma about coset graphs is well known and the proof
of the lemma follows from the definition of coset graphs.

\begin{lemma}\label{lemma-coset}
Let $\Ga=\Cos(G, H, g)$ be a coset graph. Then $\Ga$ is $G$-arc-transitive and
\begin{itemize}
\item[(1)] $\val\Ga=|H:H\cap H^g|$;
\item[(2)] $\Ga$ is connected if and only if $\l H,g \r=G$.
\item[(3)] If $\Aut\Ga$ has a subgroup $R$ acting regularly on the vertices of $\Cos(G,H,g)$, then $\Cos(G,H,g)\cong\Cay(R,S)$, where $S=R\cap HgH$.
\end{itemize}
Conversely, each $G$-arc-transitive graph $\Sigma$ is isomorphic to a coset graph
$\Cos(G, G_{v}, g)$ with $g$ satisfying the following condition$:$

\vskip0.1in
{\noindent\bf Condition:} $g$ is a $2$-element of $G$,
$g^2\in G_{v}$, $\l G_{v},g\r=G$
and $\val\Ga=|G_{v}:G_{v}\cap G_{v}^g|$, where $v\in V\Ga$.
\end{lemma}

Following the term in \cite{DFZ17},
the element $g$ satisfying the above condition
is called a {\it feasible element} to $G$ and $G_{\a}$.

A typical induction method for studying symmetric graphs is taking normal quotient
graphs. Let $\Ga$ be an $X$-vertex-transitive graph, where $X\le\Aut\Ga$.
Suppose that $X$ has a normal subgroup $N$ which is intransitive on $V\Ga$.
Denote $V_N$ the set of $N$-orbits in $V\Ga$. The {\em normal quotient graph} $\Ga_N$ defined as the
graph with vertex set $V_N$ and two $N$-orbits $B,C\in V_N$ are adjacent in $\Ga_N$ if and only if some
vertex of $B$ is adjacent in $\Ga$ to some vertex of $C$.
By \cite[Theorem 9]{Lorimer}, we have the following lemma.
\begin{lemma}\label{lor}
Let $\Ga$ be an arc-transitive graph of prime valency $p>2$ and let
$X$ be an arc-transitive subgroup of $\Aut\Ga$.
If a normal subgroup $N$ of $X$ has more than two orbits
on $V\Ga$, then $\Ga_N$ is an $X/N$-arc-transitive graph of valency
$p$ and $N$ is semiregular on $V\Ga$.Moreover, $X_v\cong(X/N)_B$ for any $v\in V\Ga$ and $B\in V\Ga_N$.
\end{lemma}

Let $\Ga$ be a graph and let $s$ be a positive integer. Recall that the graph $\Ga$ is said to be {\em $(G,s)$-arc-transitive},
if $G$ acts transitively on the set of $s$-arcs of $\Ga$, where an {\em $s$-arc} is an $(s+1)$-tuple $(v_0, v_1,\cdots,v_s)$
of $s+1$ vertices satisfying $(v_{i-1},v_i)\in E\Ga$ and $v_{i-1}\not=v_{i+1}$ for all $i$.
The graph $\Ga$ is called {\em $(G,s)$-transitive} if it is $(G,s)$-arc-transitive
but not $(G,s+1)$-arc-transitive. In particular, an $(\Aut\Ga,s)$-arc-transitive
or $(\Aut\Ga,s)$-transitive graph is just called {\em $s$-arc-transitive}
or {\em $s$-transitive graph}.
The following lemma is about the stabilizers of 13-valent symmetric graphs,
refer to \cite[Theorem 2.1]{GHL15} and \cite[Corollary 1.3]{LLL18}.

\begin{lemma}\label{stabilizer}
Let $\Ga$ be an 13-valent $(G,s)$-transitive graph, where $G\le\Aut\Ga$ and $s\ge 1$.
Let $\a\in V\Ga$. If $G_{\a}$ is soluble, then $|G_{\a}|\div 1872$.
Further, the couple $(s,G_{\a})$ lies in the following table.

\[\begin{array}{c|c} \hline
s & 1  \\ \hline
G_{\a} & \ZZ_{13},~\F_{26},~\F_{39},~\F_{52},~\F_{78},~\F_{26}\times\ZZ_2,~\F_{39}\times\ZZ_3,~\F_{52}\times\ZZ_2, ~\F_{52}\times\ZZ_4, \\& ~\F_{78}\times\ZZ_2,~\F_{78}\times\ZZ_3,~\F_{78}\times\ZZ_6 \\ \hline
s & 2 \\ \hline
G_{\a}& \F_{156},~\F_{156}\times\ZZ_2,~\F_{156}\times\ZZ_3,~\F_{156}\times\ZZ_4,~\F_{156}\times\ZZ_6 \\ \hline
s & 3 \\ \hline
G_{\a}& \F_{156}\times\ZZ_{12}  \\ \hline
\end{array}\]

If $G_{\a}$ is insoluble, then either $G_{\a}\cong \A_{13}$, $\S_{13}$, $\A_{13}\times \A_{12}$, $(\A_{13}\times \A_{12}):\ZZ_2$ or $\S_{13}\times \S_{12}$, or one of the following holds.
\begin{itemize}
	\item[(1)] $s=2$, $G_{\a}\cong ((9:\ZZ_l)\times \PSL(3,3))$, where $\ZZ_l\leq \ZZ_{2}$.
	\item[(2)] $s=2$, $G_{\a}\cong \bm{\O}_3(G_\a).\ZZ_l.\PSL(3,3)$, where $\ZZ_l\leq \ZZ_{2}$.
	\item[(3)] $s=3$, $G_{\a}\cong ((\ZZ_{3}:\ZZ_l.\PSL(2,3).O)\times \PSL(3,3))$, where $\ZZ_l\leq \ZZ_{2}$ and $O\leq \ZZ_{2}$.
	\item[(4)] $s=3$, $G_{\a}\cong \bm{\O}_3(G_\a).\ZZ_l.((\PSL(2,3).O)\times \PSL(3,3))$, where $\ZZ_l\leq \ZZ_{2}$ and $O\leq \ZZ_{2}$.
\end{itemize}
\end{lemma}

The following lemma is about primitive permutation groups of degree less than 1872, refer to \cite{Ron05}.
\begin{lemma}\label{sg}
	Let $T$ be a primitive permutation group on $\Omega$ and let $K$ be the stabilizer of a point $w\in\Omega$.
	If $T$ is a nonabelian simple group, $K$ is soluble and $|\Omega|$ divides 1872, then the triple $(T,K,|\Omega|)$ lies in the following Table \ref{table1}.

\begin{table}[ht]
	\caption{Primitive permutation groups of degree less than 1872}\label{table1}
	\vskip-0.3in
	\[\begin{array}{l|l|l|l|l|l|l|l|l}  \hline
		T&K&|\Omega|&T&K&|\Omega|&T&K&|\Omega| \\ \hline
		\A_{13}&\S_{11}&78&\A_{39}&\A_{38}&39&\A_{18}&\A_{17}&18 \\ \hline
		\PSL(2,13)&\D_{14}&78&\A_{48}&\A_{47}& 48&\A_{117}&\A_{116}&117 \\ \hline
		\PSL(4,53)&\PSp(4,3):2&117&\A_{78}&\A_{77}& 78&\A_{104}&\A_{103}&104 \\ \hline
		\PSU(3,4)&\A_5\times \ZZ_5&208&\A_{156}&\A_{155}& 156&\A_{36}&\A_{35}&36 \\ \hline
		\M_{11}&\PSL(2,11)&12&\A_{312}&\A_{311}&312&\A_{234}&\A_{233}&234 \\ \hline
		\M_{12}&\M_{11}& 12&\A_{624}&\A_{623}&624&\A_{208}&\A_{207}&208 \\ \hline
		\M_{12}&\PSL(2,11)&144&\A_{936}&\A_{935}&936&\A_{72}&\A_{71}&72 \\ \hline
		\M_{12}:2&\PSL(2,11):2&144&\A_{1872}&\A_{1871}&1872&\A_{468}&\A_{467}&468\\ \hline
		\A_{13}&\A_{12}&13&\A_{12}&\A_{11}&12&\A_{144}&\A_{143}&144  \\ \hline
		\A_{16}&\A_{15}&16&\A_{26}&\A_{25}&26 &\A_{52}&\A_{51}&52\\ \hline
		\A_{24}&\A_{23}&24 \\ \hline
	\end{array}\]
\end{table}
\end{lemma}
Let $G$ is a finite group. If $G'=G$ then $G$ is called a $perfect$  $group$, and a extension $G=N.H$ is called a $central$ $extension$ if $N\subseteq Z(G)$, the center if $G$. And $G$ is called a $covering$ $group$ of $T$ if $G$ is a perfect group and $G/Z(G)$ is isomorphic to a simple group $T$. Every nonabelian simple group $T$  has a maximal covering group, it implies that every covering group of $T$ is a factor group of the maximal covering group. The center  of the maximal covering group $G$ is the $Schur$ $multtiplier$ of $T$, denoted by $\Mult(T)$. The following lemma is about subgroups of $\ZZ_2.\A_n$, refer to \cite[Proposition 2.6]{DFZ17}.

\begin{lemma}\label{Mult}
For $n\ge 7$, all subgroups of index  $n$ in $\ZZ_2.\A_n$ are isomorphic to $\ZZ_2.\A_{n-1}$.
\end{lemma}

\begin{lemma}\label{RG}
	Let $\Ga$ be a connected $X$-arc transitive graph of valency thirteen, and $X\leq A=Aut\Ga$. Let $G\leq X$ is a regular non-abelian simple group on $V\Ga$ and let $R\not= 1$ be the soluble radical of $\A$, the largest soluble normal subgroup of $\A$. Then
 if $B=RG\not=R\times G$, then $G\lesssim GL(l,p)$ which p is a prime, integer $l\ge 2$ and $p^l\div |R|$;	
\end{lemma}
\pf Since $R$ is a solvable normal subgroup and $G$ is a non-abelian simple subgroup of $A$,we have $R\cap G \unlhd G$. It implies $R\cap G=1$, and $|B|=|R||G|$. $R$ is solvable, so $B$ has a range of normal subgroup $R_i$ such than $1=R_0<R_1<\cdot \cdot \cdot<R_s=R<B$, where $R_i\unlhd B$ and $R_{i+1}/R_i$ is ableian for $0\leq i\leq s-1$.
We assume $B=RG\not=R\times G$. Then there exists some $0\leq j\leq s-1$ so that $GR_i=G\times R_i$ for every $ {0\leq i\leq j}$, but $GR_{j+1}\not=G\times R_{j+1}$. In particular $GR_j=G\times R_j$. Since $R_j$ is solvable, $R_j\cap G=1$ and $GR_j/R_j\cong G/R_j\cap G=G$. Because $G$ is simple, we have $GR_j/R_j \cap R_{j+1}/R_j=1$, and conjugation action of $GR_j/R_j$ on $R_{j+1}/R_j$ is either trivial or faithful. Suppose the action is trivial. Then  $GR_{j+1}/R_j=(GR_j/R_j)(R_{j+1}/R_j)=GR_j/R_j\times (R_{j+1}/R_j)$, we have $GR_j\unlhd GR_{j+1}$. Noting than $G$ is characteristic in $GR_j$ as $GR_j=G\times R_j$, so $G\unlhd GR_{j+1}$, then $GR_{j+1}=G\times R_{j+1}$ which is a contradiction. It follows that this action is faithful. Since $R_{i+1}/R_i\cong \ZZ^l_p$ for some prime $p$ and integer $l$, we have $G\lesssim \GL(l,p)$ by $N/C$ theorem. And since $G$ is a non-abelian simple group, we have $l\ge2$. Obviously, it can be obtained $p^l\div |R|$. This completes the proof. \qed

\section{The proof of Theorem \ref{thm1} }
Let $\Ga=\Cay(G,S)$ be an 13-valent symmetric Cayley graph, where $G$ is a finite nonabelian simple group.
Let $A=Aut\Ga$ and let $A_{v}$ be the stabilizer of $v$ in $\A$ where $v\in V\Ga$.
Let $R$ be the soluble radical of $A$, the largest soluble normal subgroup of $A$.
Clearly, $R$ is a characteristic subgroup of $A$.
Assume that $A_v$ is soluble. Then by Lemma \ref{stabilizer}, $|A_v|$ divides 1872.

The following lemma consider the case where $R = 1$.
\begin{lemma}\label{lemma1}
Assume that $R=1$. Then $G$ is either normal in $A$ or $A$ contains a proper nonabelian simple group $T$, and $(T,G)=(\A_{13},\A_{12})$, $(\A_{39},\A_{38})$, $(\A_{117},\A_{116})$, $(\A_{208},\A_{207})$, $(\A_{312},\A_{311})$, $(\A_{936},\A_{935})$ or $ (\A_{1872},\A_{1871})$.
\end{lemma}
\pf Let $N$ be a minimal normal subgroup of $
A$. Then $N=T^d$, where $d\ge1$ and $T$ is a nonabelian simple group.
Assume that $G$ is not normal in $A$. Then since $N\cap G\unlhd G$ and $G$ is a nonabelian simple group,
$N\cap G=1$ or $G$. Assume $N\cap G=1$. Then since $A=GA_v$, we have $NG\leq A$, $|NG|=|N||G|\div |A|=|G||A_v|$, so $|N|\div|A_v|$.
It follows that $|N|\div 1872$ because $|A_v|\div 1872 $. Since $N$ is insoluble, $N$ has three divisors, by checking the simple $K_3$ groups (see \cite{Her68}), which is a contradiction. Hence $N\cap G=G$, and so $G\le N$.
If $G=N$, then $G\unlhd A$, a contradiction to the assumption. Thus $G<N$.
Assume that $d\ge2$. Then $N=T_1\times T_2\times\ldots\times T_d$ where $d\ge2$ and $T_i\cong T$ is a nonabelian simple group.
Note that $T_1\cap G\unlhd N\cap G=G$. So $T_1\cap G=1$ or $G$, if $T_1\cap G=1$, a similar argument as above, we have $|T_1|\div 1872$, which is a contradiction. Then  $T_1\cap G=G, G\leq T_1$, $|T_2|\div |N:T_1|\div |N:G|$. And $|N:G|\div |A:G|=|A_v|$,
it implies that $|T_2|\div 1872$, which is also a contradiction.
Thus, $d=1$ and $N=T$ is a nonabelian simple group.
Then $T=GT_v$, $T_v \not= 1$. Since $\Ga$ is connected and $T=N\unlhd A$,
we have $1\not=T_v^{\Ga(v)}\unlhd A_v^{\Ga(v)}$. Since $\Ga$ is $A$-arc-transitive of valency 13,
it implies that $A_v^{\Ga(v)}$ is primitive on $\Ga(v)$ and so $T_v^{\Ga(v)}$ is transitive on $\Ga(v)$ and
$13\div|T_v|$, $\Ga$ is $T$-acr-transitive of valence 13. So $|T_v|$ divides 1872.
Since $T$ has the proper subgroup $G$ with index dividing 1872, we can take a maximal proper subgroup
$K$ of $T$ which contains $G$ as a subgroup. Let $\Omega=[T : K]$. Then $|\Omega|$ divides 1872 and $T$ has
a primitive permutation representation on $\Omega$, of degree $n=|\Omega|$. Since $T$ is simple, this
representation is faithful and thus $T$ is a primitive permutation group of degree $n$. Due to the maximality of $K$ , so $K$ is the stabilizer of a point $w\in \Omega$, that is, $K = T_w$.
Consequently, by Lemma \ref{sg}, we have that the triple $(T,K,|\Omega|)$ is listed in Table \ref{table1}.
Since $|T_v|=|T:G|=|T:K||K:G|=|\Omega||K:G|$ and $|\Omega|\div1872$,
by checking the triples listed in Table \ref{table1}, we have $13$ divides $|\Omega|$.
Hence, $(T,K,|\Omega|)\not=(\M_{11},\PSL(2,11),12)$, $(\M_{12},\M_{11},12)$,$(\M_{12},\PSL(2,11),144)$, $(\M_{12}:2,\PSL(2,11):2,144)$, $(\A_{16},\A_{15},16)$,$(\A_{24},\A_{23},24)$, $(\A_{48},\A_{47},48)$,$(\A_{12},\A_{11},12)$,$(\A_{18},\A_{17},18)$,$(\A_{72},\\ \A_{71},72)$,
$(\A_{36},\A_{35},36)$ or $(\A_{144},\A_{143},144)$.

Assume that $(T,K,|\Omega|)=(\A_{13},\S_{11},78)$. Then since $G\le K$ and $G$ is a nonabelian simple group,
we have that $G$ is a proper subgroup of $K$. Since $|T:G|\div1872$ and $|\Omega|=78$, we have $|K:G|$ divides $24$. By querying the maximal subgroups of $S_{11}$, we have $G=\A_{11}$ and $|T_v|=156$. By Lemma \ref{stabilizer}, $T_v\cong\F_{156}$.By [Atlas], $T_v$ is in $\PSL(3,3)$ the maximal subgroups of $A_{13}$.
However, $\PSL(3,3)$ has no subgroup of order 156,
 a contradiction.

Assume that $(T,K,|\Omega|)=(\PSL(2,13),\D_{14},78)$. Then $K=\D_{14}$ has no simple subgroup, which is a contradiction.

Assume that $(T,K,|\Omega|)=(\PSL(4,3),\PSp(4,3):2,117)$. Then $|K:G|$ divides $16$.
By[Atlas] we have the minimum index of group $K=\PSp(4,3)$ is $27$, which is also a contradiction.

Assume that $(T,K,|\Omega|)=(\PSU(3,4),\A_{5}\times \ZZ_5,208)$. Then $|K:G|$ divides $9$, and $|K:G|=1,3,9$. Since $G$ is nonabelian simple group, no such $G$ exists, which is a contradiction.

Assume that $(T,K,|\Omega|)=(\A_{78},\A_{77}.78)$. Then $|K:G|$ divides $24$, $G=K=\A_{77}$ and $|T_v|=78$. By Lemma \ref{stabilizer}, $T_v\cong\F_{78}$. Note that $T$ has a factorization $T=GT_v$ with $G\cap T_v=G_v=1$.
By considering the right multiplication action of $T$ on the right cosets of $G$ in $T$,
we may view $T$ as a subgroup of the symmetric group $\S_n$ with $n=|T:G|=78$,
which contains a regular subgroup $T_v$.
However, $\A_{78}$ has no regular subgroup isomorphic to $\F_{78}$, a contradiction.
A similar argument, we can exclude the case $(T,K,|\Omega|)=(\A_{156},\A_{155},156)$,$(\A_{26},\A_{25},26)$, $(\A_{52},\A_{51},52)$,
$(\A_{234},\A_{233},234)$ or $(\A_{468},\A_{467},468)$.

Assume that $(T,K,|\Omega|)=(\A_{104},\A_{103},104)$ or $(\A_{624},\A_{623},624)$.
Then $G=\A_{103}$ or $\A_{623}$. By Lemma \ref{stabilizer}, $T_v=\F_{52}\times\ZZ_2$ or $\F_{156}\times\ZZ_4$.
Since $\Ga$ is $T$-arc-transitive, by Lemma \ref{lemma-coset}, we have $\Ga\cong\Cos(T,T_v,g)$ for some feasible element $g\in T$.
A direct computation by Magma \cite{Magma} shows that there is no feasible element to $T$ and $T_v$, a contradiction.

Thus, we have $(T,K)=(\A_{39},\A_{38})$, $(\A_{117},\A_{116})$, $(\A_{208},\A_{207})$ or $(\A_{13},\A_{12})$ .
For all these cases, it is easy to check that $G=K$. The lemma holds.
\qed

The following lemma consider the case $R\not=1$.

\begin{lemma}\label{lemma2}
	Assume that $G$ is not normal in $A$, $R\not=1$ and $R$ has at least three orbits on $V\Ga$. Then $RG=R\times G$.
\end{lemma}
\pf Let $B=RG$. By Lemma \ref{lor}, we have $R$ is semiregular on $V(\Ga)$ and $\Ga_R$ is an $A/R$-arc-transitive graph of valency 13, $A_v\cong(A/R)_m$ for any $v\in V(\Ga)$ and $m\in V(\Ga_N)$. So $(A/R)_m$ as a stabilizer of $\Ga_R$ is solvable. Besides, we have $G\cong B/R$  is vertex-transitive on $V(\Ga_N)$ and $G=G/R\cap G \cong GR/R=B/R\leq X/R$.
Since $R$ is the radical of $A$, so the radical of $A/R$ is trivial. According to Lemma \ref{lemma1}, we have $ B/R=G\cong T\leq S/R=:soc(A/R)$. Furthermore, $(S/R,B/R)=(\A_{n},\A_{n-1})$ with $n\ge 13$ and $n \div 1872$.

If $RG\not=R\times G$, then by lemma \ref{RG}, $G\lesssim \GL(l,p)$ for some prime $p$ , integer $l\ge2$ and $p^l\div |R|$. Due to $R\cap G\unlhd G$ and $G$ is simple, if $R\cap G = G$, $G\leq R$ and $G$ is soluble which is a contradiction. We have $R\cap G = 1$ and so $|R|\div|A_v|$. It follows that $|R|\div1872$. Especially, $p=2$, $2\leq l\leq4$ or $p=3$, $l=2$. Because $\GL(2,3)$, $\GL(2,2)$ and $\GL(3,2)$ does not have a nonabelian subgroup, and we have $r=4,p=2$ and $G\lesssim \GL(4,2)$. By Atlas \cite{Wilson}, $G=\A_5$, $\A_6$, $\A_7$, $\A_8$ or $\PSL(3,2)$, since $G\cong B/R=\A_n$ for $n\ge 13$, it is a contradiction. So  $RG=R\times G$. \qed

\begin{lemma}\label{lemma3}
Assume that $R\not=1$. Then $G$ is either normal in $\A$ or $A$ contains a proper nonabelian simple group $T$, and $(T,G)=(\A_{13},\A_{12})$, $(\A_{39},\A_{38})$, $(\A_{117},\A_{116})$, $(\A_{208},\A_{207})$, $(\A_{312},\A_{311})$, $(\A_{936},\A_{935})$ or $ (\A_{1872},\A_{1871})$.
\end{lemma}
\pf Assume that $R\not=1$ and $G$ is not normal in $A$.
Since $R\cap G\unlhd G$ and $G$ is simple, we have $|R|\div|\A_v|$.
So $|R|\div1872$.

If $R$ is transitive on $V\Ga$, then$|R:R_v|=|V\Ga|=|G|$, and $|G|\div |A_v|\div 1872$.
Since  $G$ is nonabelian simple, it is a contradiction.

If $R$ has exactly two orbits on $V\Ga$, then $\Ga$ is bipartite. It follows that the stabilizer of $G$ on the biparts is a subgroup of $G$ with index 2, which is a contradiction as $G$ is a simple group.

Thus, $R$ has more than two orbits on $V\Ga$. Let $\bar A = A/R$ and let $\bar\Ga=\Ga_R$.
By Lemma \ref{lor}, $R$ is semi-regular on $V\Ga$,
 $\bar\Ga$ is $\bar A$-arc-transitive, and so $B=R\times G$ by lemma \ref{lemma2}. Then
Let $\bar N$ be a minimal normal subgroup of $\bar A$ and let $N$ be the full preimage of $\bar N$ under $ A\rightarrow A/R$.
Since $R$ is the largest soluble normal subgroup of $ A$, we have $\bar N$ is insoluble.
Thus $\bar N=T_1\times T_2\times \ldots\times T_d=T^d$, where $T$ is a nonabelian simple group and $d\ge1$.

We first show that $d=1$. Let $\bar G=GR/R$. Then $\bar G\cong G/(G\cap R)\cong G$ is a nonabelian simple group.
Since $\bar N\cap \bar G\unlhd\bar G$, we have $\bar N\cap\bar G=1$ or $\bar G$.
If $\bar N\cap\bar G=1$, then $|\bar N|$ divides 1872, which is a contradiction with the same discussion as before.
Hence $\bar G\le\bar N$. Since $\bar G$ is simple, $|\bar G|$ must divide the order of some composition factor
of $\bar N$, that is, $|\bar G|\div|T_1|$. If $d\ge2$ then $|T_2|$ divides $|\bar N : \bar G|$ which divides $|\bar A_{\bar v}|$ with $\bar v\in \bar{\Ga}$, which is not
possible since $\bar A_{\bar v}$ divides 1872 and $T_2$ is nonabelian simple.

Now we prove that $d=1$ and $\bar N$ is a nonabelian simple group. Further, if $\bar{A}$ has another minimal normal subgroup $\bar{M}$, by the similar discussion above, we have $\bar G\le\bar M$ and $\bar M$ is simple.It follows $\bar M\bar{N}=\bar M\times \bar{N}\leq\bar{A}$ and $\bar{M}\bar{G}\leq\bar{A} $,it imples $|\bar{M}|\div 1872$, which is a contradiction. So
$\bar N$ is the unique insoluble minimal normal subgroup of $\bar A$. Assume $G$ is not normal in $A$. Since $G{\sf\,char\,}B$, $B$ is not normal in $A$, hence $G\cong B/R$ is not normal in $\bar{A}$. Let $soc(\bar{A})=\bar{N}> \bar{G}\cong G=B/R$. By Lemma \ref{lemma1}, $(N/R=\bar{N},\bar{G}\cong B/R)=(\A_{13},\A_{12})$, $(\A_{39},\A_{38})$, $(\A_{117},\A_{116})$, $(\A_{208},\A_{207})$, $(\A_{312},\A_{311})$, $(\A_{936},\A_{935})$ or $ (\A_{1872},\A_{1871})$.

Let $C=C_N(R)$, then $C\unlhd N$. Since $B=R\times G<N$, $G$ is nonabelian simple group, so $G<C$. $C\cap R=Z(R)\leq Z(C)\leq C$, then $1\not=C/(C\cap R)\cong CR/R \unlhd N/R$, since $N/R\cong \bar{N}$ is simple group, so $CR=N$ and $C=(C\cap R).\bar{N}$ is a center extension. If $C\cap R<Z(C)$, then $1\not=Z(C)/(C\cap R)\unlhd C/(C\cap R)\cong CR/R=N/R=\bar{N}$. Due to the simplicity of $\bar{N}$, we have $Z(C)=C$, a contradiction. Hence $C\cap R=Z(C)$ and $C/Z(C)\cong \bar{N}$. Now since $C'\cap Z(C)\leq Z(C')$, we have $Z(C')/(C'\cap Z(C))\unlhd C'/(C'\cap Z(C))\cong C'Z(C)/Z(C)=(C/Z(C))'=\bar{N}'=\bar{N}=C/Z(C)$. Similarly, we obtain $C'\cap Z(C)=Z(C')$, $C'=Z(C').\bar{N}$ and $C=C'Z(C)$. Furthermore, $C'=(C'Z(C))'=C''$ and $C'$ is a covering group of $\bar{N}$. Hence $C'/Z(C)=Z(C')\leq \Mult(\bar{N})$.

Since $\bar{N}=A_n$ with $n\ge 13$, By \cite [Theorem5.14] {Kleidman}, $\Mult(\bar{N})\cong \ZZ_2$, thus $Z(C')=1$ or $\ZZ_2$. If $Z(C')=\ZZ_2$, we have $C'=\ZZ_2.\A_n$. Since $\bar{G}\cap Z(C')=1$, we obtain $\bar{G}Z(C')=\bar{G}\times \ZZ_2=\A_{n-1}\times\ZZ_2$ is a subgroup of $C'$ with index $n$, which is a contradiction by lemma \ref{Mult}. So we have $Z(C')=1$, then $C'\cong \bar{N}=N/R$ is a nonabelian simple group and $C'\cap R=1$. Since $G<C$, then $G=G'<C'$. Note that $|N|=|N/R||R|=|C'||R|$ and $G<C'\unlhd N$, we have $N=C' \times R \unlhd A$. Then $soc(A/R)=N/R=(C'\times R)/R\unlhd A/R$, thus $C'\times R\unlhd A$.
 Since $C' {\sf\,char\,} R\times C'\unlhd A$, we have $C'\unlhd A$. It follows that $(\bar{N}\cong C',\bar{G}\cong G)=(\A_{13},\A_{12})$, $(\A_{39},\A_{38})$, $(\A_{117},\A_{116})$, $(\A_{208},\A_{207})$, $(\A_{312},\A_{311})$, $(\A_{936},\A_{935})$ or $ (\A_{1872},\A_{1871})$,  the lemma is true by taking $C'=T$. \qed

Now, we are ready to prove Theorem \ref{thm1}.

\vskip0.1in
{\noindent\bf Proof of Theorem \ref{thm1}.} By Lemma \ref{lemma1} and Lemma \ref{lemma3},
we have that $G$ is either normal in $\A$ or
$G=\A_{12}$, $\A_{38}$, $\A_{116}$, $\A_{207}$, $\A_{311}$, $\A_{935}$ or $\A_{1871}$, .
By \cite[Theorem 1.3]{FMW11},
for each prime $p>5$, there is a connected $p$-valent non-normal
$\A_p$-arc-transitive Cayley graph of $\A_{p-1}$,
so $\Ga$ exists for the case $G=\A_{12}$; And if $G=\A_{38}$, $\A_{116}$ and $\A_{207}$, by Examples \ref{a39}, \ref{a117} and \ref{a208} below,
there exist connected
$13$-valent symmetric non-normal
Cayley graphs of $\A_{n}$ with $n=38$, $116$ or $207$,
the last statement of Theorem \ref{thm1} is true.
This completes the proof of Theorem \ref{thm1}.\qed
\section{The examples and the full automorphism groups}
In this section, we construct some examples to show that, for $G=\A_{38}$, $\A_{116}$ or $\A_{207}$,
there exist non-normal 13-valent symmetric Cayley graphs of $G$ and determine the full automorphism group of these graphs.

\begin{example}\rm \label{a39}
Let $X$ be the group consisting of all even permutations in  $\Omega_1 ={\lbrace1,2,...,39\rbrace}$ and $G$ be the group consisting of all even permutations in  $\Omega_2 ={\lbrace2,3,...,39\rbrace}$,then $X\cong \A_{39}$ and $G\cong \A_{38}$.
\begin{itemize}
\item[$x$]=(1, 2, 4)(3, 6, 11)(5, 9, 16)(7, 13, 22)(8, 14, 24)(10, 18, 29)(12, 20, 27)(15, 26, 21)(17, 28, 34)(19, 25, 35)(23, 32, 37)(30, 38, 39)(31, 36, 33),
\item[$y$]=(1, 3, 7, 14, 25, 29, 26, 36, 38, 16, 27, 37, 28)(2, 5, 10, 6, 12, 21, 13, 23, 33, 24, 34, 39, 35)(4, 8, 15, 9, 17, 22, 18, 30, 32, 11, 19, 31, 20),
\item[$g$]=(1, 7)(2, 22)(3, 5)(4, 13)(6, 16)(9, 11)(14, 24)(18, 29)(20, 27)(21, 26)(23, 31)(25, 35)(28, 34)(32, 33)(36, 37)(38, 39).
\end{itemize}

Let $H=\l x,y \r$ and let $\Ga=\Cos(X,H,g)$.

\vskip0.1in By Magma \cite{Magma}, $H=\l y\r:\l x \r\cong\F_{39}$, $\l H,g\r=X$
and $|H:H \cap H^g|=13$.
By Lemma \ref{lemma-coset}(1)(2), $\Ga$ is a connected $\A_{39}$-arc-transitive $13$-valent graph.
Also, it is easy to see that $H$ is regular on $\lbrace1,2,...,39\rbrace$.
Hence the vertex stabilizer $X_1=G\cong\A_{38}$ is regular on $V\Ga=[X:H]$,
that is, $\Ga$ is a Cayley graph of $\A_{38}$.
Finally, since $G\cong\A_{38}$ is not normal in $X\cong\A_{39}$,
we have that $\Ga$ is non-normal.
\end{example}

\begin{example}\rm \label{a117}
Let $X$ be the group consisting of all even permutations in  $\Omega_1 ={\lbrace1,2,...,117\rbrace}$ and $G$ be the group consisting of all even permutations in  $\Omega_2 ={\lbrace2,3,...,117\rbrace}$,then $X\cong \A_{117}$ and $G\cong \A_{116}$.
\begin{itemize}
\item[$x$]=(1, 2, 3)(4, 10, 16)(5, 11, 17)(6, 12, 18)(7, 19, 22)(8, 20, 23)(9, 21, 24)(13, 28, 31)(14, 29, 32)(15, 30, 33)(25, 37, 34)(26, 38, 35)(27, 39, 36)(40, 41, 42)(43, 49, 55)(44, 50, 56)(45, 51, 57)(46, 58, 61)(47, 59, 62)(48, 60, 63)(52, 67, 70)(53, 68, 71)(54, 69, 72)(64, 76, 73)(65, 77, 74)(66, 78, 75)(79, 80, 81)(82, 88, 94)(83, 89, 95)(84, 90, 96)(85, 97, 100)(86, 98, 101)(87, 99, 102)(91, 106, 109)(92, 107, 110)(93, 108, 111)(103, 115, 112)(104, 116, 113)(105, 117, 114);
\item[$y$]=(1, 71, 99, 13, 78, 85, 23, 56, 82, 30, 64, 90, 35, 40, 110, 21, 52, 117, 7, 62, 95, 4, 69, 103, 12, 74, 79, 32, 60, 91, 39, 46, 101, 17, 43, 108, 25, 51, 113)(2, 72, 97, 14, 76, 86, 24, 57, 83, 28, 65, 88, 36, 41, 111, 19, 53, 115, 8, 63, 96, 5, 67, 104, 10, 75, 80, 33, 58, 92, 37, 47, 102, 18, 44, 106, 26, 49, 114)(3, 70, 98, 15, 77, 87, 22, 55, 84, 29, 66, 89, 34, 42, 109, 20, 54, 116, 9, 61, 94, 6, 68, 105, 11, 73, 81, 31, 59, 93, 38, 48, 100, 16, 45, 107, 27, 50, 112);
\item[$g$]=(2, 40)(3, 79)(4, 55)(10, 94)(11, 44)(12, 45)(17, 83)(18, 84)(19, 46)(20, 47)(21, 48)(22, 85)(23, 86)(24, 87)(26, 27)(28, 52)(29, 53)(30, 54)(31, 91)(32, 92)(33, 93)(34, 103)(35, 105)(36, 104)(37, 64)(38, 66)(39, 65)(42, 80)(49, 82)(56, 89)(57, 90)(61, 97)(62, 98)(63, 99)(70, 106)(71, 107)(72, 108)(73, 115)(74, 117)(75, 116)(77, 78)(113, 114)(1, 41)(2, 42)(3, 40)(4, 94)(5, 50)(6, 51)(7, 58)(8, 59)(9, 60)(10, 82)(11, 56)(12, 57)(13, 67)(14, 68)(15, 69)(16, 88)(17, 44)(18, 45)(19, 61)(20, 62)(21, 63)(22, 46)(23, 47)(24, 48)(25, 76)(26, 78)(27, 77)(28, 70)(29, 71)(30, 72)(31, 52)(32, 53)(33, 54)(34, 64)(35, 66)(36, 65)(37, 73)(38, 75)(39, 74)(104, 105)(113, 114)(116, 117).
\end{itemize}

Let $H=\l x,y\r$ and let $\Ga=\Cos(X,H,g)$.

\vskip0.1in By Magma \cite{Magma}, $H\cong \ZZ_3 \times \F_{39}$, $\l H,g\r=X$
and $|H:H \cap H^g|=13$.
Hence Lemma \ref{lemma-coset} implies that $\Ga$ is a connected $\A_{117}$-arc-transitive $13$-valent graph.
Also, with a similar discussion as above,
we have that $H$ is regular on $\lbrace1,2,...,117\rbrace$,
and $\Ga$ is a non-normal Cayley graph of $G=\A_{116}$.
\end{example}

\begin{example}\rm \label{a208}
	Let $X$ be the group consisting of all even permutations in  $\Omega_1 ={\lbrace1,2,...,208\rbrace}$ and $G$ be the group consisting of all even permutations in  $\Omega_2 ={\lbrace2,3,...,208\rbrace}$,then $X\cong \A_{208}$ and $G\cong \A_{207}$.
\begin{itemize}
\item[$x$]=(1, 2, 3, 4)(5, 13, 11, 21)(6, 14, 12, 22)(7, 15, 9, 23)(8, 16, 10, 24)(17, 37, 27, 41)(18, 38, 28, 42)(19, 39, 25, 43)(20, 40, 26, 44)(29, 34, 49, 46)(30, 35, 50, 47)(31, 36, 51, 48)(32, 33, 52, 45)(53, 54, 55, 56)(57, 65, 63, 73)(58, 66, 64, 74)(59, 67, 61, 75)(60, 68, 62, 76)(69, 89, 79, 93)(70, 90, 80, 94)(71, 91, 77, 95)(72, 92, 78, 96)(81, 86, 101, 98)(82, 87, 102, 99)(83, 88, 103, 100)(84, 85, 104, 97)(105, 106, 107, 108)(109, 117, 115, 125)(110, 118, 116, 126)(111, 119, 113, 127)(112, 120, 114, 128)(121, 141, 131, 145)(122, 142, 132, 146)(123, 143, 129, 147)(124, 144, 130, 148)(133, 138, 153, 150)(134, 139, 154, 151)(135, 140, 155, 152)(136, 137, 156, 149)(157, 158, 159, 160)(161, 169, 167, 177)(162, 170, 168, 178)(163, 171, 165, 179)(164, 172, 166, 180)(173, 193, 183, 197)(174, 194, 184, 198)(175, 195, 181, 199)(176, 196, 182, 200)(185, 190, 205, 202)(186, 191, 206, 203)(187, 192, 207, 204)(188, 189, 208, 201);
\item[$y$]=(1, 77, 136, 192, 22, 92, 109, 165, 42, 68, 150, 206, 17, 53, 129, 188, 36, 74, 144, 161, 9, 94, 120, 202, 50, 69, 105, 181, 32, 88, 126, 196, 5, 61, 146, 172, 46, 102, 121, 157, 25, 84, 140, 178, 40, 57, 113, 198, 16, 98, 154, 173)(2, 78, 133, 189, 23, 89, 110, 166, 43, 65, 151, 207, 18, 54, 130, 185, 33, 75, 141, 162, 10, 95, 117, 203, 51, 70, 106, 182, 29, 85, 127, 193, 6, 62, 147, 169, 47, 103, 122, 158, 26, 81, 137, 179, 37, 58, 114, 199, 13, 99, 155, 174)(3, 79, 134, 190, 24, 90, 111, 167, 44, 66, 152, 208, 19, 55, 131, 186, 34, 76, 142, 163, 11, 96, 118, 204, 52, 71, 107, 183, 30, 86, 128, 194, 7, 63, 148, 170, 48, 104, 123, 159, 27, 82, 138, 180, 38, 59, 115, 200, 14, 100, 156, 175)(4, 80, 135, 191, 21, 91, 112, 168, 41, 67, 149, 205, 20, 56, 132, 187, 35, 73, 143, 164, 12, 93, 119, 201, 49, 72, 108, 184, 31, 87, 125, 195, 8, 64, 145, 171, 45, 101, 124, 160, 28, 83, 139, 177, 39, 60, 116, 197, 15, 97, 153, 176);
\item[$g$]=(1, 54)(3, 158)(4, 106)(5, 66)(6, 65)(7, 67)(8, 68)(9, 171)(10, 172)(11, 170)(12, 169)(13, 14)(17, 89)(18, 90)(19, 91)(20, 92)(21, 118)(22, 117)(23, 119)(24, 120)(25, 195)(26, 196)(27, 193)(28, 194)(29, 86)(30, 87)(31, 88)(32, 85)(41, 141)(42, 142)(43, 143)(44, 144)(45, 137)(46, 138)(47, 139)(48, 140)(49, 190)(50, 191)(51, 192)(52, 189)(55, 157)(56, 105)(57, 58)(61, 163)(62, 164)(63, 162)(64, 161)(73, 110)(74, 109)(75, 111)(76, 112)(77, 175)(78, 176)(79, 173)(80, 174)(93, 121)(94, 122)(95, 123)(96, 124)(97, 136)(98, 133)(99, 134)(100, 135)(101, 185)(102, 186)(103, 187)(104, 188)(107, 160)(113, 179)(114, 180)(115, 178)(116, 177)(125, 126)(129, 199)(130, 200)(131, 197)(132, 198)(153, 202)(154, 203)(155, 204)(156, 201)(167, 168).
\end{itemize}
	
Let $H=\l x,y\r$ and let $\Ga=\Cos(X,H,g)$.

\vskip0.1in By Magma \cite{Magma}, $H\cong \ZZ_4 \times \F_{52}$, $\l H,g\r=X$
and $|H:H \cap H^g|=13$.
Hence Lemma \ref{lemma-coset} implies that $\Ga$ is a connected $\A_{208}$-arc-transitive $13$-valent graph.
Also, with a similar discussion as above,
we have that $H$ is regular on $\lbrace1,2,...,208\rbrace$,
and $\Ga$ is a non-normal Cayley graph of $G=\A_{207}$.
\end{example}

At the end of this paper, we determine the full automorphism group
of the graph constructed in Example \ref{a39}.
Recall that a transitive permutation group is called {\it quasiprimitive}
if each of its minimal normal subgroups is transitive.

\begin{lemma}\label{rem-d10-s10}
Let $\Ga=\Cos(X,H,g)$ be as in Example~$\ref{a39}$. Then $\Aut\Ga\cong\A_{39}$ or $\S_{39}$
and $\Ga$ is $1$-transitive.
\end{lemma}

\pf Recall that $\A_{38}\cong G<X\cong\A_{39}$
and $\Ga$ is a connected $X$-arc-transitive $13$-valent Cayley graph of $G$.
Let $\A=\Aut\Ga$ and $v\in V\Ga$.
By \cite[Theorem 2.1]{GHL15} and \cite[Corollary 1.3]{LLL18}, $|\A_{v}|\div2^{20}\cdot3^{10}\cdot5^4\cdot7^2\cdot11^2\cdot13$.

Assume $\A$ is not quasiprimitive on $V\Ga$.
Then $\A$ has an intransitive minimal normal subgroup $N$.
Set $F=NX$.
Since $X$ is nonabelian simple and $N\cap X\lhd X$, we have $N\cap X=1$ or $X$.
If $N\cap X=X$, then $N$ is transitive on $V\Ga$, a contradiction.
Suppose $N\cap X=1$. Then $F=N:X$ and $|N|=|F:X|$ divides $|\A:X|$. Since
$|V\Ga|=|\A:\A_{v}|=|X:X_{v}|$, we have
$|\A:X|=|\A_v:X_v|$ divides $2^{20}\cdot3^9\cdot5^4\cdot7^2\cdot11^2$,
so is $|N|$.
Since $|V\Ga|=|G|=|\A_{38}|$, if $N$ has exactly two orbits on $V\Ga$. It follows that the stabilizer of $G$ on the biparts is a subgroup of $G$ with index 2, which is a contradiction as $G$ is a simple group. So $N$ has at least three orbits on $V\Ga$.
By Lemma \ref{lor}, $N$ is semi-regular on $V\Ga$,
and so $|N|$ divides $|V\Ga|=|\A_{38}|$.

Suppose that $N$ is insoluble. Note that $|N|\div 2^{20}\cdot3^9\cdot5^4\cdot7^2\cdot11^2$. Then by checking the simple $K_3$ groups (see \cite{Her68}), the simple $K_4$ groups (see \cite[Theorem 1]{BCM01}) and the simple $K_5$ groups (see \cite[Theorem A]{A. Jafarzadeh})
we can conclude that $N\cong\A_5$, $\A_5^2$, $\A_5^3$, $\A_5^4$, $\A_6$, $\A_6^2$, $\A_6^3$, $\A_6^4$, $\A_7$, $\A_7^2$, $\A_8$, $\A_8^2$, $\A_9$, $\A_9^2$, $\A_{10}$, $\A_{10}^2$, $\A_{11}$, $\A_{11}^2$, $\A_{12}$, $\PSL(2,7)$, $\PSL(2,7)^2$,  $\PSL(2,8)$, $\PSL(2,8)^2$, $\PSL(2,11)$, $\PSL(2,49)$, $\PSU(3,3)$, $\PSU(3,3)^2$, $\PSL(3,4)$, $\PSL(3,4)^2$,
$\PSU(3,5)$, $\PSU(4,2)$, $\PSU(4,2)^2$, $\PSU(4,3)$, $\PSU(5,2)$, $\PSU(6,2)$, $\PSp(6,2)$, $\PSp(6,2)^2$, $\PSO(7,2)$, $\PSO^{+}(8,2)$, $\PSO(7,2)^2$, $\M_{11}$, $\M_{11}^2$, $\M_{12}$, $\M_{12}^2$, $\M_{22}$, $\M_{22}^2$, $\J_2$, $\J_2^2$, $\HS$,$\McL$.
Then since $|N||\A_{39}|=|N||X|=|F|=|V\Ga||F_{v}|=|\A_{38}||F_{v}|$,
we have $|F_{v}|=39\cdot|N|$.
By checking the orders of the stabilizers of connected 13-valent symmetric graphs given in Lemma \ref{stabilizer},
none of these values for $|F_v|$ satisfies the orders, a contradiction.

Now suppose that $N$ is soluble. Noting that $|N|\div |\A_v:X_{v}|$,
$|\A_v:X_{v}|\div 2^{20}\cdot3^9\cdot5^4\cdot7^2\cdot11^2$,
we have $N\cong\ZZ_2^i$, $\ZZ_3^j$, $\ZZ_5^k$, $\ZZ_7^m$ or $\ZZ_{11}^n$, where $1\le i\le 20$, $1\le
j\le9$, $1\le k\le4$, $1\le m\le2$ and $1\le n\le2$.
Note that $\N_F(N)/\C_F(N)=F/\C_F(N)\lesssim \Aut(N)\cong \GL(i,2)$, $\GL(j,3)$, $\GL(k,5)$, $\GL(m,7)$ or $\GL(n,11)$.
Clearly, $N\le \C_F(N)$. If $N=\C_F(N)$, then $\A_{39}\cong X\cong F/N=F/\C_F(N)
\lesssim \GL(i,2)$, $\GL(j,3)$, $\GL(k,5)$, $\GL(m,7)$ or $\GL(n,11)$. However, by
Magma \cite{Magma}, each of $\GL(i,2)$, $\GL(j,3)$, $\GL(k,5)$, $\GL(m,7)$ and $\GL(n,11)$ has no subgroup
isomorphic to $\A_{39}$ for $1\le i\le 20$,$1\le
j\le9$, $1\le k\le4$, $1\le m\le2$ and $1\le n\le2$, a contradiction.
Hence $N<\C_F(N)$ and $1\not=\C_F(N)/N\unlhd F/N\cong\A_{39}$. It
follows $F=\C_F(N)=N\times X$,
$F_v/X_v\cong F/X\cong N$,
and $F_v$ is soluble because $X_{v}\cong\F_{39}$.
By Lemma \ref{stabilizer}, we conclude that
$F_v\cong \ZZ_2\times \F_{39}$, $\ZZ_2^2\times\F_{39}$ or  $\ZZ_3{\times}\F_{39}$.
A direct computation by Magma \cite{Magma} shows that
there is no feasible element to $F$ and $F_{v}$,
it is also a contradiction.

Thus, $\A$ is quasiprimitive on $V\Ga$.
Let $M$ be a minimal normal subgroup of $\A$. Then $M=T^d$,
with $T$ a nonabelian simple group, is transitive on
$V\Ga$, so $|V\Ga|=|\A_{38}|$ divides $|M|$ and $37\div|T|$.
If $d\ge2$, then $37^2\div|M|$,
which is a contradiction because $|\A|\div|\A_{38}|\cdot2^{20}\cdot3^{10}\cdot5^4\cdot7^2\cdot11^2\cdot13$
is not divisible by $37^2$.
Hence $d=1$ and $M=T\lhd \A$. Let $C=\C_\A(T)$. Then $C\lhd \A$
and $CT = C{\times}T$. If $C\not=1$, then $C$ is transitive on $V\Ga$ as
$\A$ is quasiprimitive on $V\Ga$,
with a similar discussion as above, we have $C$ is insoluble and $37\div|C|$.
Therefore, $37^2\div|CT|\div|A|$, again a contradiction.
Hence $C=1$ and $\A$ is almost simple.

Since $M\cap X\unlhd X\cong\A_{39}$,
we have $M\cap X=1$ or $X$. If $M\cap X=1$, then
$|M|\div 2^{20}\cdot3^9\cdot5^4\cdot7^2\cdot11^2$, it is a contradiction as $|\A_{38}|\div|M|$.
Thus, $M\cap X=X$ and so $\A_{39}\cong X\le M$. Hence $M$ is a nonabelian simple
group satisfying $|\A_{39}|\div|M|$ and $|M|\div|\A_{38}|\cdot2^{20}\cdot3^{10}\cdot5^4\cdot7^2\cdot11^2\cdot13$.
By \cite[P.135--136]{Gorenstein}, we can conclude that $M=X\cong\A_{39}$.
Thus $\A\le\Aut(\M)\cong\S_{39}$. If $\A\cong\S_{39}$, then $|\A_v|=|\A:G|=78$,
and so $\A_v\cong\F_{78}$ by Lemma \ref{stabilizer}.
A direct computation by Magma \cite{Magma} shows that
there is feasible element to $\A$ and $\A_{v}$. Hence $\A\cong\S_{39}$ or $\A_{39}$
and $\Ga$ is $1$-transitive. \qed

\end{document}